\documentclass[11pt]{article}
\usepackage[a4paper,margin=26mm]{geometry}
\usepackage{amsmath,amssymb,amsthm}
\usepackage{mathtools}
\usepackage{microtype}
\usepackage{enumitem}
\usepackage[mathlines]{lineno}
\usepackage{graphicx,float}
\usepackage{caption}
\usepackage{xurl}
\usepackage[hidelinks]{hyperref}

\newtheorem{theorem}{Theorem}[section]
\newtheorem{lemma}[theorem]{Lemma}
\newtheorem{proposition}[theorem]{Proposition}

\theoremstyle{remark}
\newtheorem{remark}[theorem]{Remark}

\newcommand{\vol}{\operatorname{vol}}
\newcommand{\diam}{\operatorname{diam}}
\newcommand{\conv}{\operatorname{conv}}

\title{\bfseries The Longest-Edge Bisection Algorithm May Produce Degenerating Simplices in Higher Dimensions}
\author{Sergey Korotov$^{1}$ \qquad Michal K\v{r}\'i\v{z}ek$^{2}$\\[0.75em]
\small $^{1}$ Department of Business and Mathematics, M\"alardalen University,\\
\small Box 883, 721 23 V\"aster\aa s, Sweden\\[0.3em]
\small $^{2}$ Institute of Mathematics, Czech Academy of Sciences,\\
\small \v{Z}itn\'a 25, CZ-115 67 Prague 1, Czech Republic\\[0.3em]
\small Corresponding author: Sergey Korotov, \href{mailto:sergey.korotov@mdu.se}{sergey.korotov@mdu.se}\\
\small Michal K\v{r}\'i\v{z}ek: \href{mailto:krizek@math.cas.cz}{krizek@math.cas.cz}}
\date{}

\begin{document}
\maketitle

\begin{abstract}
For every integer $d\ge3$ and $a\in(0,1)$, we give a one-parameter infinite family $E_d(a)$ of
$d$-simplices for which exactly $d-1$ classical longest-edge bisections produce a simplex
congruent to $E_d(a/2)$. Starting from the two-step tetrahedral recurrence, we
obtain the family of $d$-simplices simply by adjoining $d-3$ equal transverse vertices and
performing one additional bisection for each added dimension. Thus
$d-1=2+(d-3)$: an active tetrahedral core is followed by transverse completion. The proof uses
only explicit coordinates, midpoint formulas, and squared-distance comparisons.
The normalized $d$-volume tends to zero, and the same branch violates both the
$d$-dimensional minimum- and maximum-angle conditions. Higher-dimensional degeneration
itself is already known; the contribution here is a short and fully explicit construction requiring
no orbit or projective-dynamical machinery. A short appendix embeds the local branch into
conforming simplicial partitions with vanishing maximal diameter.
Consequently, classical longest-edge bisection does not guarantee uniform shape regularity
of simplicial partitions in dimensions three and higher. This is particularly relevant to the
geometric assumptions commonly imposed in the convergence analysis of finite element methods.
\end{abstract}


\section{Introduction}
Under classical longest-edge bisection (LEB), at each step a longest edge $PQ$ of a $d$-simplex is selected,
its midpoint $M=(P+Q)/2$ is inserted, and one of the two resulting children is selected.
When several edges attain the maximum length, the classical LEB rule permits an arbitrary
choice among them. For $d=2,3,\ldots$ and a nondegenerate $d$-simplex $T$, its shape regularity is measured by the following ratio:
\[
  \vartheta_d(T):=\frac{\vol_d(T)}{(\diam T)^d}.
\]
The same scale-invariant quantity is used by Adiprasito, Kalmanovich, and
Solomon~\cite{AKS}: an LEB branch degenerates when $\vartheta_d$ tends to zero as the refinement proceeds. From the
finite element viewpoint, a uniform positive lower bound for $\vartheta_d$ is a standard
volume regularity condition. In arbitrary dimension it is equivalent to the generalized
Zl\'amal angle condition; see Brandts, Korotov, and K\v{r}\'i\v{z}ek~\cite{BKK2011}.
For triangles and tetrahedra, equivalence with the usual angle and ball conditions was
established in~\cite{BKK2008}.

Although LEB guarantees diameter decay in arbitrary space dimension~\cite{Kearfott}, the
question of whether it preserves shape regularity in dimension three remained open for a long time.
Recent results on the orbits of tetrahedra generated by LEB are presented in~\cite{MichaudKorotov}.
A counterexample based on a two-step recurrence was constructed in~\cite{Korotov2026}, yielding
an infinite sequence of tetrahedra $T_n$ whose normalized volumes tend to zero, i.e.\ $\vartheta_3(T_n)\to0$ as $n\to\infty$.

The same two-step tetrahedral construction subsequently appeared in~\cite[Section~3.2]{AKS},
where its identity with the example in~\cite{Korotov2026} is explicitly acknowledged.
Adiprasito, Kalmanovich, and Solomon also established further results using a dynamical/projective
framework. In dimension four they give a strictly LEB-admissible nine-bisection return and
prove an open-set result; a face-lifting lemma then yields higher-dimensional degeneration.

The purpose of the present note is different: we give a substantially simpler, fully explicit
higher-dimensional construction. Starting from an elementary tetrahedral counterexample given in~\cite{Korotov2026}, we
simply adjoin $d-3$ equal transverse vertices, one for each additional dimension. Accordingly,
the number of bisections decomposes as
\[
  d-1=2+(d-3).
\]
The first two bisections reproduce the tetrahedral recurrence, while one further bisection per
added dimension restores the prescribed configuration of the transverse vertices. Ties among
the longest edges have no significance for the present counterexample, since classical LEB
permits an arbitrary choice among them and it suffices to construct one admissible degenerating branch.

\begin{remark}[Why the lifting mechanism starts in dimension three]\label{rem:2d-lift}
The face-lifting lemma of Adiprasito, Kalmanovich, and Solomon~\cite[Lemma~3.9]{AKS}
allows a degenerating LEB branch of a $k$-dimensional face of a $d$-simplex to be lifted to a
degenerating branch of the entire simplex, where $2\le k<d$. It therefore does not fail for $k=2$, $d=3$;
rather, its hypothesis is unavailable. Classical planar LEB has a uniform lower angle
bound~\cite{RosenbergStenger}, so no planar LEB branch can have $\vartheta_2\to0$.
Hence there is no degenerating two-dimensional branch to use as a seed. Dimension three is
the first dimension in which such a seed exists.
\end{remark}

A conforming-partition realization of the local branch is recorded in
Appendix~\ref{app:conforming}. In particular, diameter decay alone does not ensure that the
geometric hypotheses of the usual shape-regular finite element error estimates hold.

\section{Construction and main result}
\begin{figure}[H]
\centering
\includegraphics{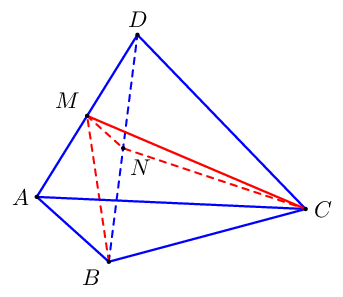}
\caption{The tetrahedron $ABCD$, its child $MBCD$, and its grandchild $MBCN$, where $M$ and $N$ are the midpoints of $AD$ and $BD$, respectively. The red edges indicate the two successive longest-edge bisections}
\label{fig:core}
\end{figure}

Let $d\ge3$ and let $e_1,\dots,e_d$ be the canonical basis of $\mathbb R^d$. Fix
\begin{equation}\label{eq:a-range}
  0<a<1,
\end{equation}
and set
\begin{equation}\label{eq:vertices}
\begin{aligned}
  A&=ae_3, & B&=0,\\
  C&=\sqrt a\,e_1, &
  D&=\frac{\sqrt a}{2}e_1+\frac{\sqrt{7a}}{2}e_2,\\
  N&=\frac{B+D}{2}=\frac D2, &
  X_j&=N+a e_j \qquad(4\le j\le d).
\end{aligned}
\end{equation}
For $d=3$ there are no additional vertices $X_j$. Define
\begin{equation}\label{eq:Ed}
  E_d(a)=\conv\{A,B,C,D,X_4,\dots,X_d\}.
\end{equation}
Using $B$ as the base vertex and subtracting $N=D/2$ from each column $X_k$
($4\le k\le d$) of the edge-vector matrix, we obtain
\begin{equation}\label{eq:volume}
  \vol_d(E_d(a))=\frac{\sqrt7}{2d!}\,a^{d-1}>0.
\end{equation}
The exponent $d-1$ reflects the anisotropic scaling of the construction: the tetrahedron $ABCD$
contributes the factor $a^2$, while each of the $d-3$ transverse directions contributes
one further factor $a$. Hence $a^2a^{d-3}=a^{d-1}$.
Put
\begin{equation}\label{eq:MZ}
  M=\frac{A+D}{2},\qquad Z_j=\frac{C+X_j}{2}\qquad(4\le j\le d).
\end{equation}
The first two bisections of the tetrahedral core are illustrated in Fig.~\ref{fig:core}.

The next theorem will be proved in Section~3.

\begin{theorem}[$(d-1)$-step recurrence]\label{thm:recurrence}
Let $d\ge3$ and $0<a<1$. For any permutation $\sigma$ of $\{4,\dots,d\}$, consider the directed
word
\begin{equation}\label{eq:word}
  AD,\quad BD,\quad CX_{\sigma(4)},\quad CX_{\sigma(5)},\quad\dots,\quad CX_{\sigma(d)}.
\end{equation}
For $d=3$, the transverse part of \eqref{eq:word} is absent, so the word is simply
$AD,BD$. In all dimensions the displayed word is LEB-admissible. At the first step select the child containing $D$; at the second select the
child containing $B$; at every transverse step select the child containing $C$. The final
simplex is
\begin{equation}\label{eq:return}
  \conv\{M,N,B,C,Z_4,\dots,Z_d\}\cong E_d(a/2).
\end{equation}
Initially $AD$ is uniquely longest. After the first bisection, $BD$ and $CD$ are exactly the
longest edges and are tied. After the second bisection, and throughout the transverse
completion as long as unprocessed vertices remain, the edges $CX_j$ corresponding to the as-yet unprocessed transverse vertices are
exactly the longest edges and are mutually tied. Thus the transverse bisections may be
performed in any order. Repetition yields an infinite degenerating branch in every dimension
$d\ge3$.
\end{theorem}

The return word displays the mechanism directly:
\[
  \underbrace{AD,\ BD}_{\text{two-step tetrahedral core}}
  \quad\text{and}\quad
  \underbrace{CX_{\sigma(4)},\dots,CX_{\sigma(d)}}_{\text{transverse completion}}.
\]

\begin{remark}[Tie breaking is optional]\label{rem:tiebreak}
If one wants unique transverse choices, replace $X_j=N+a e_j$ by
$X_j=N+a\lambda_j e_j$ with
\[
  \frac12<\lambda_d<\cdots<\lambda_4<1.
\]
The same squared-distance argument gives the ordered return
$AD,BD,CX_4,\dots,CX_d$ with all transverse selected edges uniquely longest. This
multi-parameter refinement is not needed for the counterexample and will not be used below.
\end{remark}

\section{Proof of the recurrence}
\subsection{The two-step tetrahedral core}
\begin{lemma}\label{lem:firsttwo}
The initial edge $AD$ is uniquely longest. After replacing $A$ by $M$, the longest edges are
exactly $BD$ and $CD$, each of squared length $2a$.
\end{lemma}

\begin{proof}
The six core squared-edge lengths, in the order $(AB,AC,AD,BC,BD,CD)$, are
\begin{equation}\label{eq:corelengths}
  (a^2,\ a+a^2,\ 2a+a^2,\ a,\ 2a,\ 2a).
\end{equation}
For every added vertex $X_j$,
\begin{equation}\label{eq:Xlengths}
\begin{aligned}
 |AX_j|^2&=\frac a2+2a^2,
 & |BX_j|^2=|DX_j|^2&=\frac a2+a^2,\\
 |CX_j|^2&=a+a^2,
 & |X_iX_j|^2&=2a^2\qquad(i\ne j).
\end{aligned}
\end{equation}
Since $0<a<1$, every quantity in \eqref{eq:Xlengths} is smaller than
$2a+a^2=|AD|^2$. Thus $AD$ is uniquely longest.

After bisecting $AD$, the selected simplex has vertices $M,B,C,D,X_4,\dots,X_d$.
Among the core vertices,
\[
  |BD|^2=|CD|^2=2a,
\]
while
\[
  |MB|^2=|MD|^2=\frac a2+\frac{a^2}{4},\qquad
  |MC|^2=a+\frac{a^2}{4}.
\]
The unchanged edges $BC$, $BX_j$, and $DX_j$ are also shorter than $BD$ by
\eqref{eq:corelengths} and \eqref{eq:Xlengths}. Further,
\[
  |MX_j|^2=\frac54a^2,
  \qquad |CX_j|^2=a+a^2,
  \qquad |X_iX_j|^2=2a^2.
\]
All these are smaller than $2a$ because $0<a<1$. Hence the only longest edges are
$BD$ and $CD$. Choosing $BD$, bisecting it at $N$, and selecting the child containing $B$ leaves the simplex
with vertices
\[
  M,N,B,C,X_4,\dots,X_d.
\]
\end{proof}

\subsection{Transverse completion}
\begin{lemma}\label{lem:transverse}
Starting from $\{M,N,B,C,X_4,\dots,X_d\}$, as long as unprocessed transverse vertices
remain, the edges $CX_j$ belonging to them are exactly the longest edges. Hence they may be bisected in
any order, always selecting the child containing $C$.
\end{lemma}

\begin{proof}
Let $I$ and $J$ be the processed and unprocessed index sets, respectively. The current simplex has core
$K=\{M,N,B,C\}$, vertices $Z_i$ for $i\in I$, and vertices $X_j$ for $j\in J$.
If $J=\varnothing$, the transverse completion is finished. Otherwise, for arbitrary $j\in J$ put
\[
  L:=|CX_j|^2=a+a^2.
\]
Within the core,
\begin{equation}\label{eq:Kmax}
  \max_{U,V\in K}|UV|^2=|MC|^2=a+\frac{a^2}{4}<L.
\end{equation}
For $i,k\in I$, $i\ne k$, and $\ell,m\in J$, $\ell\ne m$, the remaining squared distances are
\begin{equation}\label{eq:completion-lengths}
\begin{aligned}
 |MX_\ell|^2&=\frac54a^2, & |NX_\ell|^2&=a^2,
 & |BX_\ell|^2&=\frac a2+a^2,\\
 |MZ_i|^2&=\frac a4+\frac{a^2}{2},
 & |NZ_i|^2&=\frac a4+\frac{a^2}{4},
 & |BZ_i|^2&=\frac a2+\frac{a^2}{4},\\
 |CZ_i|^2&=\frac a4+\frac{a^2}{4},
 & |X_\ell X_m|^2&=2a^2,
 & |X_\ell Z_i|^2&=\frac a4+\frac54a^2,\\
 |Z_iZ_k|^2&=\frac{a^2}{2}.&&&
\end{aligned}
\end{equation}
Each equality applies whenever the indices occurring in it exist.
For $0<a<1$, every quantity in \eqref{eq:Kmax}--\eqref{eq:completion-lengths} is
smaller than $L$. Thus the only longest edges are the mutually tied edges $CX_j$ with
$j\in J$. Bisecting any one of them replaces $X_j$ by
$Z_j=(C+X_j)/2$ and preserves the same description for the remaining indices. Induction
completes the transverse part.
\end{proof}

\subsection{Exact return}
\begin{proof}[Proof of Theorem~\ref{thm:recurrence}]
Lemmas~\ref{lem:firsttwo} and~\ref{lem:transverse} prove the admissibility of the prescribed LEB steps. Relabel the final
core tetrahedron by
\begin{equation}\label{eq:relabel}
  A'=M,\qquad B'=N,\qquad C'=B,\qquad D'=C.
\end{equation}
Its squared lengths, in the order $(A'B',A'C',A'D',B'C',B'D',C'D')$, are
\[
  \left(\frac{a^2}{4},\ \frac a2+\frac{a^2}{4},\ a+\frac{a^2}{4},\
  \frac a2,\ a,\ a\right),
\]
which are exactly \eqref{eq:corelengths} with $a$ replaced by $a/2$. Moreover, by \eqref{eq:MZ}, \eqref{eq:vertices}, and \eqref{eq:relabel} we get
\begin{equation}\label{eq:Zreturn}
  Z_j=\frac{C+X_j}{2}
     =\frac{C+N}{2}+\frac a2 e_j
     =\frac{B'+D'}{2}+\frac a2 e_j.
\end{equation}
The isometry between the two tetrahedral cores extends to an isometry of $\mathbb R^d$
whose orthogonal part fixes $e_4,\ldots,e_d$. Thus the transverse vertices also have exactly
the defining form of $E_d(a/2)$. This proves \eqref{eq:return}.
\end{proof}

\section{Degeneration and angle conditions}
\begin{theorem}[Quantified degeneration]\label{thm:degeneration}
For $0<a<1$,
\begin{equation}\label{eq:diameter}
  (\diam E_d(a))^2=a(2+a),
\end{equation}
and
\begin{equation}\label{eq:theta}
  \vartheta_d(E_d(a))
  =\frac{\sqrt7}{2d!}\frac{a^{(d-2)/2}}{(2+a)^{d/2}}
  \xrightarrow[a\downarrow0]{}0.
\end{equation}
The same degeneration holds along all intermediate simplices of the infinite branch. For
$d=3$, \eqref{eq:volume}, \eqref{eq:diameter}, and \eqref{eq:theta} reduce exactly to the
corresponding tetrahedral formulas in~\cite{Korotov2026}.
\end{theorem}

\begin{proof}
The volume formula is \eqref{eq:volume}, and Lemma~\ref{lem:firsttwo} shows that $AD$ is
the initial longest edge, which gives \eqref{eq:diameter}. If $S_r$ is the intermediate simplex
after $r$ bisections of one recurrence cycle, $0\le r\le d-1$, then
\[
  \vol_d(S_r)=2^{-r}\vol_d(E_d(a)).
\]
Every $S_r$ contains the vertices $B$ and $C$, whose distance is $\sqrt a$, while $S_r\subset E_d(a)$. Hence,
\[
  \sqrt a\le\diam S_r\le\sqrt{a(2+a)}.
\]
Therefore,
\[
  0<\vartheta_d(S_r)
  \le 2^{-r}\frac{\sqrt7}{2d!}\,a^{(d-2)/2}
  \le \frac{\sqrt7}{2d!}\,a^{(d-2)/2}
  \xrightarrow[a\downarrow0]{}0,
\]
which proves degeneration at all positions of the cycle.
\end{proof}

\begin{remark}[Step count and dimension]\label{rem:steps}
From \eqref{eq:volume} we have
\[
  \frac{\vol_d(E_d(a/2))}{\vol_d(E_d(a))}=2^{-(d-1)}.
\]
Let $a_0\in(0,1)$ be given. Since every bisection halves the $d$-volume, exactly $d-1$ bisections are forced for this
return. Thus higher dimensions require more individual bisections to complete one exact
contraction $a\mapsto a/2$: after $k$ full cycles, $a_k=2^{-k}a_0$ and $k(d-1)$
bisections have occurred. The degenerating mechanism itself, however, is already present in
the tetrahedral core; the additional $d-3$ steps only complete the transverse geometry.
\end{remark}

\begin{proposition}[Explicit minimum- and maximum-angle degeneration]\label{prop:angles}
At the recurrence-cycle endpoints, the triangular subsimplex $ABC$ has
\begin{equation}\label{eq:alpha-min}
  \alpha(a):=\angle ACB=\arctan\sqrt a,
  \qquad \sin\alpha(a)=\sqrt{\frac{a}{1+a}},
\end{equation}
so $\alpha(a)\to0$. After the first bisection, the current $d$-simplex contains the
tetrahedral subsimplex
\[
  O(a)=\conv\{M,B,C,D\},
\]
whose interior dihedral angle $\Theta(a)$ at $MC$ satisfies
\begin{equation}\label{eq:Theta-max}
  \cos\Theta(a)=\frac{2a-7}{\sqrt{(4a+7)(8a+7)}}\longrightarrow-1.
\end{equation}
Hence, $\Theta(a)\to\pi$. Therefore, the branch violates both the $d$-dimensional minimum-
and maximum-angle conditions.
\end{proposition}

\begin{proof}
From \eqref{eq:vertices},
\[
  \overrightarrow{CA}=(-\sqrt a,0,a,0,\ldots,0),\qquad
  \overrightarrow{CB}=(-\sqrt a,0,0,\ldots,0),
\]
which gives \eqref{eq:alpha-min}. Since $ABC$ is a vertex-generated subsimplex, this
violates the minimum-angle condition in the sense of Korotov and K\v{r}\'i\v{z}ek
\cite{KorotovKrizek2025}; see also~\cite{BKK2011,KorotovVatne2020}.

The first bisection replaces only $A$ by $M$, so $O(a)$ occurs as a vertex-generated
tetrahedral subsimplex in every dimension. Formula \eqref{eq:Theta-max} is exactly the
dihedral-angle formula obtained for the three-dimensional counterexample in~\cite{Korotov2026}.
Hence, the maximum-angle condition of~\cite{KorotovKrizek2025} fails; equivalently, the
generalized maximum-angle condition of Khademi, Korotov, and Vatne
\cite[Theorem~3.7]{KhademiKorotovVatne2019} is not valid.
\end{proof}

\begin{remark}\label{rem:three-failures}
Along the same branch,
\[
  \vartheta_d(E_d(a))\to0,\qquad \alpha(a)\to0,\qquad \Theta(a)\to\pi,
\]
with
\[
  \alpha(a)\sim\sqrt a,\qquad
  \pi-\Theta(a)\sim\frac4{\sqrt7}\sqrt a.
\]
Thus volume degeneration and both angle-condition failures are visible through explicit
one-parameter formulas.
\end{remark}

\section{Four dimensions and relation to previous work}
For $d=4$ there is a single transverse vertex,
\[
  X_4=N+a e_4
      =\left(\frac{\sqrt a}{4},\frac{\sqrt{7a}}{4},0,a\right).
\]
The exact return is
\[
  E_4(a)\xrightarrow{\,AD\,}\cdot
  \xrightarrow{\,BD\,}\cdot
  \xrightarrow{\,CX_4\,}E_4(a/2),
\]
with only the core tie $BD=CD$, and
\[
  \vol_4(E_4(a))=\frac{\sqrt7}{48}a^3.
\]

The tetrahedral core is exactly the one-parameter squared-edge family appearing in
\cite{Korotov2026}. After division by $a$, its squared-edge vector, in the order
$(AD,AB,AC,DB,DC,BC)$, is
\[
  (2+a,\ a,\ 1+a,\ 2,\ 2,\ 1),
\]
which is also the family in~\cite[Eq.~(17)]{AKS}, up to relabeling.

The two-step tetrahedral counterexample underlying the present construction was introduced
in~\cite{Korotov2026}, submitted to arXiv on 24 August 2026. The identical construction
appears in~\cite[Section~3.2]{AKS}, submitted to arXiv on 8 September 2026, where this
coincidence is explicitly acknowledged. The authors of~\cite{AKS} also establish further
results, including tie-free degeneration and the existence of open sets of degenerating simplices.
In particular, their four-dimensional hyperbolic construction uses a strictly LEB-admissible
nine-bisection return; their Lemma~3.9 lifts degenerating face branches to arbitrary higher
dimensions. The contribution of the present note is simplicity and explicitness: starting from
the original tetrahedral counterexample~\cite{Korotov2026}, an explicit degenerating
$d$-simplex is obtained merely by adjoining $d-3$ equal transverse vertices and performing one
additional LEB step per added dimension. The resulting exact $(d-1)$-step recurrence is verified
by direct squared-distance calculations and yields closed volume and angle formulas, with no
orbit or projective-dynamical machinery.

A different all-dimensional self-similarity was described by Moshagen
\cite[Theorem~2.1]{Moshagen2026}. For the unit simplex $\conv\{0,e_1,\dots,e_d\}$,
a particular admissible sequence of $d$ longest-edge bisections gives a homothetic descendant
with all edge lengths halved. Indeed, select successively the edges $e_1e_2,\ldots,e_1e_d$
and then $0e_1$, each time selecting the child containing $e_1$. The first $d-1$
selected edges have length $\sqrt{2}$, and the last has length $1$, which are the
respective longest-edge lengths. The resulting simplex is
homothetic to the initial one with center $e_1$ and ratio $1/2$. Repeating this sequence
preserves shape, whereas the present $(d-1)$-step return changes shape and degenerates.

To the best of our knowledge, no earlier one-parameter exact degenerating recurrence
\[
  E_d(a)\longmapsto E_d(a/2)
\]
in $d-1$ classical LEB steps, valid for every $d\ge4$, has been reported yet. No claim of first higher-dimensional
degeneration, genericity, or tie-free refinement is made.

\appendix
\section{Conforming-partition realization}\label{app:conforming}
The preceding recurrence concerns one descendant branch. The following auxiliary statement
places it inside conforming simplicial partitions with vanishing maximal diameter.

\begin{proposition}[Global-edge realization]\label{prop:partition}
Let $\Omega=E_d(a_0)$ with $0<a_0<1$. Then there exists a sequence of finite conforming (i.e. face-to-face)
simplicial partitions $\{\mathcal T_n\}_{n\ge0}$ of $\overline\Omega$ such that:
\begin{enumerate}[label=\textup{(\roman*)},leftmargin=8mm]
\item every refinement step selects an edge of globally maximal length and bisects it in every
simplex incident to it;
\item $h_n:=\max_{T\in\mathcal T_n}\diam T\to0$;
\item the partitions contain an infinite nested branch whose recurrence-cycle endpoints are
congruent to
\[
  E_d(a_0),\quad E_d(a_0/2),\quad E_d(a_0/2^2),\quad\dots.
\]
Hence the family is not uniformly shape regular.
\end{enumerate}
\end{proposition}

\begin{proof}
If a globally longest edge $PQ$ has length $c$, its midpoint $M$ satisfies, for every other
vertex $X$ of an incident simplex,
\begin{equation}\label{eq:median}
  |X-M|^2=\frac{|X-P|^2+|X-Q|^2}{2}-\frac{|P-Q|^2}{4}
  \le\frac{3c^2}{4}.
\end{equation}
Thus every new edge created by bisecting an edge of length $c$ has length at most
$\rho c$, where $\rho=\sqrt3/2<1$. Assign generation $0$ to initial edges. When an edge of
generation $q$ is bisected, assign generation $q+1$ to every newly created edge.
Each edge is bisected at most once and generates only finitely many new edges, because
all partitions are finite. Hence each generation contains only finitely many edges,
and an edge of generation $q$ has length at most $\rho^qH_0$, where $H_0$ is the
largest initial edge. Hence, for every $L>0$, only finitely many edges of length at least $L$
can occur. Repeated globally longest-edge refinement therefore drives the maximum edge
length, and hence the maximum simplex diameter, to zero.

Now suppose the current distinguished simplex has prescribed next branch edge $e$ of length
$L$. First refine globally longest edges of length larger than $L$, always in every
incident simplex. This waiting phase is finite and does not alter the distinguished simplex,
whose longest edge has length $L$. Then $e$ is globally longest; bisect it in every incident
simplex and select, inside the distinguished simplex, the child prescribed by
Theorem~\ref{thm:recurrence}. Bisecting the same edge in all incident simplices preserves
face-to-face conformity. Repetition embeds the complete recurrent branch and proves the
three assertions.
\end{proof}

\section*{Statements and Declarations}

\noindent\textbf{Funding.} The work of M. K\v{r}\'i\v{z}ek was supported by grant no. 24-10586S of the Czech Science Foundation and by the Czech Academy of Sciences (RVO 67985840).

\noindent\textbf{Competing interests.} The authors declare that they have no competing interests.

\noindent\textbf{Author contributions.} Both authors contributed equally to all aspects of this work. Both authors read and approved the final manuscript.

\noindent\textbf{Data availability.} This work is theoretical; no datasets were generated or analysed.

\noindent\textbf{Use of generative AI.} No AI tools were used. The authors are responsible for all the mathematical results and the final text of the paper.

\end{document}